\documentclass[12pt]{article}
\usepackage{amsfonts}
\usepackage{bbm}
\usepackage{mathrsfs}
\usepackage{amssymb}

\def\hang{\hangindent\parindent}
\def\textindent#1{\indent\llap{#1\enspace}\ignorespaces}

\title{Computation of Minimal Homogeneous Generating Sets \\ and Minimal Standard Bases for Ideals of Free Algebras}

\author{Huishi Li\\
{\small Department of Applied Mathematics},
{\small College of Information Science and Technology}\\
{\small Hainan University}, {\small  Haikou 570228, China}}

\date{}

\begin{document}
\maketitle
\begin{center}
\begin{minipage}{120mm}\def\KX{K\langle X\rangle}\def\LH{{\bf LH}}\def\G{{\cal G}}
{\small {\bf Abstract.} Let $\KX =K\langle X_1,\ldots ,X_n\rangle$
be the free algebra generated by $X=\{ X_1,\ldots ,X_n\}$ over a
field $K$. It is shown that with respect to any weighted
$\mathbb{N}$-gradation attached to $\KX$,  minimal homogeneous
generating sets for finitely generated graded two-sided ideals of
$\KX$ can be algorithmically computed, and that if an ungraded
two-sided ideal $I$ of $\KX$ has a finite Gr\"obner basis $\G$ with
respect to a graded monomial ordering on $\KX$, then a minimal
standard basis for $I$ can be computed via computing a minimal
homogeneous generating set of the associated graded ideal
$\langle\LH (I)\rangle$.}
\end{minipage}\end{center} {\parindent=0pt\vskip 6pt

{\bf 2010 Mathematics subject classification} Primary 16W70;
Secondary 16Z05.\vskip 6pt

{\bf Key words} Free Algebra, Homogeneous generating set, Gr\"obner
basis, Standard basis.}


\def\v5{\vskip .5truecm}\def\QED{\hfill{$\Box$}}\def\hang{\hangindent\parindent}
\def\textindent#1{\indent\llap{#1\enspace}\ignorespaces}
\def\item{\par\hang\textindent}
\def \r{\rightarrow}\def\OV#1{\overline {#1}}
\def\normalbaselines{\baselineskip 24pt\lineskip 4pt\lineskiplimit 4pt}
\def\mapdown#1{\llap{$\vcenter {\hbox {$\scriptstyle #1$}}$}
                                \Bigg\downarrow}
\def\mapdownr#1{\Bigg\downarrow\rlap{$\vcenter{\hbox
                                    {$\scriptstyle #1$}}$}}
\def\mapright#1#2{\smash{\mathop{\longrightarrow}\limits^{#1}_{#2}}}
\def\NZ{\mathbb{N}}

\def\LH{{\bf LH}}\def\LM{{\bf LM}}\def\LT{{\bf
LT}}\def\KX{K\langle X\rangle}
\def\B{{\cal B}} \def\LC{{\bf LC}} \def\G{{\cal G}} \def\FRAC#1#2{\displaystyle{\frac{#1}{#2}}}
\def\SUM^#1_#2{\displaystyle{\sum^{#1}_{#2}}} \def\O{{\cal O}}  \def\J{{\bf J}}\def\BE{\B (\mathbbm{e})}
\def\PRCVE{\prec_{\varepsilon\hbox{-}gr}}\def\S{{\cal S}}
\def\HL{{\rm LH}}\def\NB{\mathbb{N}_{\mathscr{B}}}


\def\LH{{\bf LH}}\def\LM{{\bf LM}}\def\LT{{\bf
LT}}\def\KX{K\langle X\rangle}
\def\B{{\cal B}} \def\LC{{\bf LC}} \def\G{{\cal G}} \def\FRAC#1#2{\displaystyle{\frac{#1}{#2}}}
\def\SUM^#1_#2{\displaystyle{\sum^{#1}_{#2}}} \def\O{{\cal O}}  \def\J{{\bf J}}\def\BE{\B (\mathbbm{e})}
\def\PRCVE{\prec_{\varepsilon\hbox{-}gr}}\def\S{{\cal S}}

\def\HL{{\rm LH}}\def\NB{\mathbb{N}_{\mathscr{B}}}
\vskip 1truecm

\section*{1. Introduction and Preliminary}
Throughout this paper, $K$ denotes a commutative field, $K^*=K-\{
0\}$,  algebras are meant associative $K$-algebras. Unless otherwise
stated, ideals of algebras are meant two-sided ideals, and an ideal
generated by a subset $S$ is denoted by $\langle S\rangle$.
Moreover, we use $\mathbb{N}$ to denote the additive monoid of
nonnegative integers. \v5

Let $\KX =K\langle X_1,\ldots ,X_n\rangle$ be the noncommutative
free $K$-algebra generated by $X=\{ X_1,\ldots ,X_n\}$,  and $\B=\{
1,~X_{i_1}\cdots X_{i_s}~|~X_{i_j}\in X,~s\ge 1\}$ the standard
$K$-basis of $\KX$. For convenience, elements of $\B$ are referred
to as {\it monomials} and denoted by lower case letters $w, u,v,
s,\ldots$ . Equip $\KX$ with a weighted $\NZ$-gradation $\KX
=\oplus_{q\in\NZ}\KX_q$ by assigning each $X_i$ a {\it positive
degree} $d_{\rm gr}(X_i)=m_i$, $1\le i\le n$, that is, for each
$w=X_{i_1}\cdots X_{i_s}\in\B$, $d_{\rm gr}(w)=d_{\rm
gr}(X_{i_1})+\cdots +d_{\rm gr}(X_{i_s})=m_{i_1}+\cdots +m_{i_s}$,
and for each $q\in\NZ$, $\KX$ has the degree-$q$ homogeneous part
$\KX_q=K$-span$\{ w\in\B~|~d_{\rm gr}(w)=q\}$. If $f\in\KX_q$ is a
nonzero homogeneous element of degree $q$, then we write $d_{\rm
gr}(f)=q$.\par

Since the classical Betti number defined for {\it  graded one-sided
ideals} of a noncommutative connected $\mathbb{N}$-graded
$K$-algebra has been shown to have  a {\it two-sided version} for
 graded two-sided ideals of a monoid graded local ring ([Li1],
Proposition 3.5), it follows that if $I$ is a finitely generated
{\it  graded two-sided ideal} of $\KX$, then any two minimal
homogeneous generating sets of $I$ have the same number of
generators, and any two minimal homogeneous generating sets of $I$
contain the same number of homogeneous elements of degree $n$ for
all $n\in\NZ$. Also we know that any finitely presented connected
graded $K$-algebra $A$ is isomorphic to a quotient algebra such as
$\KX /I$, i.e., the generators of $I$ give rise to the defining
relations of $A$. In this paper we first show that the methods and
algorithms, developed in  ([CDNR], [KR]) for computing minimal
homogeneous generating sets of graded submodules in free modules
over commutative polynomial algebras, can be extended to compute
minimal homogeneous generating sets of $I$ (Section 2). Secondly, in 
consideration of the relation with standard bases of ideals in 
$\KX$,  we show that if an {\it ungraded ideal} $I$ of $\KX$ has a 
finite Gr\"obner basis $\G$ with respect to a graded monomial 
ordering $\prec_{gr}$, then a minimal standard basis of $I$, which 
has similar properties as a minimal homogeneous generating set does, 
can  be computed via computing a minimal homogeneous generating set 
of the associated graded ideal $\langle \LH (I)\rangle$ of $I$ 
(Section 3).  \v5

Concerning the Gr\"obner basis theory for the free $K$-algebra
$\KX=K\langle X_1,\ldots ,X_n\rangle$, we now recall from ([Mor],
[Gr]) some basic facts as follows. Let $\prec$ be a monomial
ordering on $\B$, which is by definition a well-ordering $\prec$ on
$\B$ satisfying: $u\prec v$ implies $wus\prec wvs$ for all
$w,u,v,s\in\B$; $v\ne u$ and $v=wus$ implies $u\prec v$ for all
$u,v,w,s\in\B$. If $f\in \KX$ is such that
$f=\sum_{i=1}^m\lambda_iw_i$ with $\lambda_i\in K^*$, $w_i\in\B$,
and $w_1\prec w_2\prec\cdots\prec w_m$, then we write $\LM (f)=w_m$
for the leading monomial of $f$, and we write $\LC (f)=\lambda_m$
for the leading coefficient of $f$. \par

Let $u,v\in \B$. we say that $u$ {\it divides} $v$, denoted $u~|~v$,
if $v=wus$ for some $w,s\in\B$. As in the commutative case, if a
monomial ordering $\prec$ on $\B$ is given, then the division of
monomials extends to a division algorithm of dividing an element $f$
by a finite subset of nonzero elements $G=\{g_1,\ldots ,g_t\}$ in
$\KX$, which gives rise to a representation
$f=\sum_{i,j}\lambda_{ij}w_{ij}g_iu_{ij}+r$, where $\lambda_{ij}\in
K$, $w_{ij},u_{ij}\in \B$, $g_i\in G$, satisfying $\LM
(w_{ij}g_iu_{ij})\preceq\LM (f)$ for all $\lambda_{ij}\ne 0$, and if
$r\ne 0$ such that $r=\sum_k\mu_kv_k$ with $\mu_k\in K^*$, $v_k\in
\B$, then $\LM (r)\preceq\LM (f)$ and $\LM (g_j){\not |}~v_k$ for
all $k$. We write $\OV f^G=r$ and call it a {\it remainder} of $f$
on division by $G$. If $\OV f^G=0$, then we say that $f$ is {\it
reduced to zero} on division by $G$. A nonzero element $f\in \KX$ is
said to be {\it normal} (mod $G$) if $f=\OV f^G$. Moreover, a subset
$G$ of nonzero elements in $\KX$ is said to be LM-{\it reduced} if
$\LM (g_i){\not |}~\LM (g_j)$ for all $g_i\ne g_j$ in $G$. \par

Given a monomial ordering $\prec$ on $\B$ and a subset $\G$ of
nonzero elements in $\KX$, let $I=\langle\G\rangle$ be the ideal of
$\KX$ generated by $\G$. If for any nonzero element $f\in I$, there
is a $g_i\in\G$ such that $\LM (g_i)|\LM (f)$, then $\G$ is called a
{\it Gr\"obner basis} of $I$. For a graded ideal $I$ of $\KX$, a
Gr\"obner basis $\G$ of $I$ consisting of homogeneous elements is
called a {\it homogeneous Gr\"obner basis} of $I$. A Gr\"obner basis
$\G$ is said to be {\it minimal} if $\LM (g_i){\not |}~\LM (g_j)$
for all $g_i\ne g_j$ in $\G$.\par

Let $f,g\in\KX$ be two nonzero elements. If there are monomials
$u,v\in\B$ such that {\parindent=.75truecm\par
\item{(1)} $\LM (f)u=v\LM (g)$, and
\item{(2)} $\LM (f){\not |}~v$ and $\LM (g)\not |~u$,}{\parindent=0pt\par

then the element
$$o(f,u;~v,g)=\frac{1}{\LC (f)}(f\cdot u)-
\frac{1}{\LC (g)}(v\cdot g)$$ is referred to as an overlap element
of $f$ and $g$.}\v5

The next theorem and the following algorithm are known as the
implementation of Bergman's diamond lemma [Ber].{\parindent=0pt\v5

{\bf Theorem} (Termination theorem in the sense of  [Gr]) Let $\G
=\{ g_1,\ldots ,g_m\}$ be an LM-reduced subset of $\KX$. then $\G$
is a Gr\"obner basis for the ideal $I=\langle\G\rangle$ if and only
if for each pair $g_i,g_j\in\G$, including $g_i=g_j$, every overlap
element $o(g_i,u;~v,g_j)$ of $g_i$ and $g_j$ has the property
$\overline{o(g_i,u;~v,g_j)}^{\G}=0,$ that is, $o(g_i,u;~v,g_j)$ is
reduced to 0 by the division by $\G$.\par\QED}\par

If a given LM-reduced subset $G=\{g_1,\ldots ,g_t\}$ of $\KX$ is not
a Gr\"obner basis for the ideal $I=\langle G\rangle$, then the very
noncommutative version of the Buchberger Algorithm (cf. [Mor], [Gr])
computes a (possibly infinite) Gr\"obner basis for $I$. For the use
of next section we recall this algorithm as
follows.{\parindent=0pt\par

\underline{{\bf Algorithm 1}
~~~~~~~~~~~~~~~~~~~~~~~~~~~~~~~~~~~~~~~~~~~~~~~~~~~~~~~~~~
~~~~~~~~~~~~~~~~~~~~~~~~~~~~~~~~~~~~~~~~~~~~~}{\par

\textsc{INPUT:} $G_0=\{ g_1,...,g_t\}$\par

\textsc{OUTPUT:} $\G =\{ g_1,...,g_m,...\}$, a Gr\"obner basis for
$I$\par

\textsc{INITIALIZATION:} $\G :=G_0$, $\O :=\{ 
o(g_i,g_j)~|~g_i,g_j\in G_0\}$\par

\textsc{BEGIN}\par

~~~~\textsc{WHILE} $\O\ne\emptyset$ \textsc{DO}\par

~~~~~~~~Choose any $o(g_i,g_j)\in \O$\par

~~~~~~~~$\O :=\O-\{o(g_i,g_j)\}$\par

~~~~~~~~$\overline{o(g_i,g_j)}^{\G}=r$\par

~~~~~~~~\textsc{IF} $r\ne 0$ \textsc{THEN}\par

~~~~~~~~~~~~$\O :=\O\cup \{o(g,r), ~o(r,g),~o(r,r)~|~g\in\G\}$ \par

~~~~~~~~~~~~$\G :=\G\cup\{ r\}$}\par

~~~~~~~~\textsc{END}{\parindent=0pt\par

~~~~$\vdots$

\underline{~~~~~~~~~~~~~~~~~~~~~~~~~~~~~~~~~~~~~~~~~~~~~~~~~~~~~~~~~~
~~~~~~~~~~~~~~~~~~~~~~~~~~~~~~~~~~~~~~~~~~~~~~~~~~~~~~~~~~~~~~~~}}}\v5

\section*{2. Computation of Minimal Homogeneous Generating Sets}
Let $\KX =K\langle X_1,\ldots ,X_n\rangle$ be the free $K$-algebra
generated by $X=\{ X_1,\ldots ,X_n\}$ and $\B$ the standard
$K$-basis of $\KX$. Fix a weighted $\NZ$-gradation $\KX
=\oplus_{q\in\NZ}\KX_q$ for $\KX$ by assigning  each $X_i$ a
positive degree $d_{\rm gr}(X_i)=m_i$, $1\le i\le n$. Let $\prec$ be
a monomial ordering on $\B$. Based on {\bf Algorithm 1} presented in
Section 1, 1n this section we show that  the methods and algorithms,
developed in ([CDNR], [KR]) for computing minimal homogeneous
generating sets of graded submodules in free modules over
commutative polynomial algebras, can be adapted for computing
minimal homogeneous generating sets of a finitely generated graded
two-sided  ideal $I$ of $\KX$.  All notions, notations, and
conventions given in Section 1 are maintained. {\parindent=0pt\v5

{\bf 2.1. Definition} Let $G=\{ g_1,\ldots ,g_t\}$ be a subset of
homogeneous elements of $\KX$, $I=\langle G\rangle$ the graded ideal
generated by $G$, and let $n\in\NZ$, $G_{\le n}=\{ g_j\in G~|~d_{\rm
gr}(g_j)\le n\}$. If, for each nonzero homogeneous element $f\in I$
with $d_{\rm gr}(f)\le n$, there is some $g_i\in G_{\le n}$ such
that $\LM (g_i)|\LM (f)$ with respect to $\prec$, then we call
$G_{\le n}$ an $n$-{\it truncated Gr\"obner basis} of $I$.}\v5

Noticing that every $w\in\B$ is a homogeneous element of $\KX$,
verification of the  lemma below is straightforward.
{\parindent=0pt\v5

{\bf 2.2. Lemma} Let $\G =\{ g_1,\ldots ,g_t\}$ be a homogeneous
Gr\"obner basis for the graded ideal $I=\langle \G\rangle$ of $\KX$
with respect to the given monomial ordering $\prec$ on $\B$. For
each $n\in\NZ$, put $\G_{\le n}=\{ g_j\in\G~|~d_{\rm gr}(g_j)\le
n\}$, $I_{\le n}=\cup_{q=0}^nI_q$ where each $I_q$ is the degree-$q$
homogeneous part of $I$, and let $I(n)=\langle I_{\le n}\rangle$ be
the graded ideal generated by $I_{\le n}$.  The following statements
hold.\par

(i) $\G_{\le n}$ is an $n$-truncated Gr\"obner basis of $I$. Thus,
if $f\in\KX$ is a homogeneous element with $d_{\rm gr}(f)\le n$,
then $f\in I$ if and only if $\OV{f}^{\G_{\le n}}=0$, i.e., $f$ is
reduced to zero on division by $\G_{\le n}$.\par

(ii) $I(n)=\langle\G_{\le n}\rangle$, and $\G_{\le n}$ is an
$n$-truncated Gr\"obner basis of $I(n)$.\par\QED\v5

{\bf Convention} In what follows, we let $o(f,g)$ represent {\it any
overlap element} of two nonzero elements $f,g\in\KX$.}\v5

In light of {\bf Algorithm 1}, an $n$-truncated Gr\"obner basis is
characterized as follows.{\parindent=0pt\v5

{\bf 2.3. Proposition} Let $I=\langle G\rangle$ be the graded ideal
of $\KX$ generated by a finite set of nonzero homogeneous elements
$G=\{ g_1,\ldots ,g_m\}$. Without loss of generality, we assume that
$G$ is LM-reduced (see Section 1). For each $n\in\NZ$, put $G_{\le
n}=\{ g_j\in G~|~d_{\rm gr}(g_j)\le n\}$. The following statements
are equivalent with respect to the given monomial ordering $\prec$
on $\B$.\par

(i) $G_{\le n}$ is an $n$-truncated Gr\"obner basis of $I$.\par

(ii) For each $(g_i,g_j)\in G\times G$, every overlap element
$o(g_i,g_j)$ of $d_{\rm gr}(o(g_i,g_j))\le n$ is reduced to zero on
division by $G_{\le n}$, i.e., $\OV{o(g_i,g_j)}^{G_{\le
n}}=0$.\vskip 6pt

{\bf Proof} Recall that if $\LM (g_i)=vw$ and $\LM (g_j)=wu$ for
some $u,v,w\in\B$ with $w\ne 1$, then the corresponding overlap
element of $g_i$ and $g_j$ is
$$o(g_i,u;~v,g_j)=\frac{1}{\LC (g_i)}g_iu-\frac{1}{\LC (g_j)}vg_j$$
which is obviously a homogeneous element in $I$. If $d_{\rm
gr}(o(g_i,g_j))\le n$, then it follows from (i) that (ii) holds.}
\par

Conversely, suppose that (ii) holds. To see that $G_{\le n}$ is an
$n$-truncated Gr\"obner basis of $I$, let us run ({\bf Algorithm 1})
with the initial input data $G$. Without optimizing {\bf Algorithm
1} we may certainly assume that $G\subseteq\G$, thereby $G_{\le
n}\subseteq\G_{\le n}$, where $\G$ is the new input set returned
after a certain pass through the WHILE loop. On the other hand, by
the construction of $o(g_i,g_j)$  we know that if $d_{\rm
gr}(o(g_i,g_j))\le n$, then $d_{\rm gr}(g_i)\le n$, $d_{\rm
gr}(g_j)\le n$. Hence, the assumption (ii) implies that {\bf
Algorithm 1} does not give rise to any new element of degree $\le n$
for $\G$. Therefore, $G_{\le n}=\G_{\le n}$. By Lemma 2.2 we
conclude that $G_{\le n}$ is an $n$-truncated Gr\"obner basis of
$I$.\QED {\parindent=0pt\v5

{\bf 2.4. Corollary} Let  $I=\langle G\rangle$ be the graded ideal
of $\KX$ generated by a finite set of nonzero homogeneous elements
$G=\{ g_1,\ldots ,g_m\}$. Suppose that $G_{\le n}=\{ g_j\in
G~|~d_{\rm gr}(g_j)\le n\}$ is an $n$-truncated Gr\"obner basis of
$I$ with respect to the given monomial ordering $\prec$ on $\B$.
\par

(i) If $g\in \KX$ is a nonzero homogeneous element of $d_{\rm
gr}(g)=n$ such that $\LM (g_i){\not |}~\LM (g)$ for all $g_i\in
G_{\le n}$, then $G'=G_{\le n}\cup\{ g\}$ is an $n$-truncated
Gr\"obner basis for both the graded ideals $I'=I+\langle g\rangle$
and $I''=\langle G'\rangle$ of $\KX$.\par

(ii) If $n\le n_1$ and $g\in \KX$ is a nonzero homogeneous element
of $d_{\rm gr}(g)=n_1$ such that $\LM (g_i){\not |}~\LM (g)$ for all
$g_i\in G_{\le n}$, then $G'=G_{\le n}\cup\{ g\}$ is an
$n_1$-truncated Gr\"obner basis for the graded ideal $I'=\langle
G'\rangle$ of $\KX$.\vskip 6pt

{\bf Proof} If $g\in\KX$ is a nonzero homogeneous element of $d_{\rm
gr}(g)=n_1\ge n$ and $\LM (g_i){\not |}~\LM (g)$ for all $g_i\in
G_{\le n}$, then it is straightforward to see that $d_{\rm gr}(H)>n$
for every nonzero $H\in\{ o(g_i,g),o(g,g_i),o(g,g)~|~g_i\in G\}$.
Hence both (i) and (ii) hold by Proposition 2.3. \QED\v5

{\bf 2.5. Proposition} (Compare with ([KR], Proposition 4.5.10)) 
Given a finite set of nonzero homogeneous elements $F =\{ f_1,\ldots 
,f_m\}\subset \KX$ with $d_{\rm gr}(f_1)\le d_{\rm gr}(f_2)\le 
\cdots\le d_{\rm gr}(f_m)$, and a positive integer $n_0\ge d_{\rm 
gr}(f_1)$, the following algorithm computes an $n_0$-truncated 
Gr\"obner basis $\G =\{ g_1,\ldots ,g_t\}$ for the graded ideal 
$I=\langle F\rangle$ of $\KX$, such that $d_{\rm gr}(g_1)\le d_{\rm 
gr}(g_2)\le\cdots\le d_{\rm gr}(g_t)$. {\parindent=0pt\newpage

\underline{\bf Algorithm 2
~~~~~~~~~~~~~~~~~~~~~~~~~~~~~~~~~~~~~~~~~~~~~~~~~~~~~~~~~~~~~~~~~~~~~~~~~~~~~~~~~~~~~~~~~~}\par

$\begin{array}{l} \textsc{INPUT}:~ F= \{ 
f_1,...,f_m\}~\hbox{with}~d_{\rm gr}(f_1)\le d_{\rm gr}(f_2)\le 
\cdots\le d_{\rm gr}(f_m);~ n_0~\hbox{with}~n_0\ge d_{\rm 
gr}(f_1)\\
\textsc{OUTPUT}: ~\G =\{
g_1,...,g_t\},~\hbox{an}~n_0\hbox{-truncatedGr\"obner basis of}~I\\
\textsc{INITIALIZATION}:~ \O_{\le n_0} :=\emptyset ,~W :=F,~\G :=\emptyset ,~t' :=0\\
\textsc{LOOP}\\
n:=\min\{d_{\rm gr}(f_i),~d_{\rm gr}(o(g_{\ell},g_q))~|~f_i\in W
,~o(g_{\ell},g_q)\in\O_{\le n_0}\}\\
\O_n:=\{ o(g_{\ell},g_q)\in\O_{\le n_0}~|~d_{\rm gr}(o(g_{\ell},g_q))=n\},~W_n:=\{ f_j\in W~|~d_{\rm gr}(f_j)=n\}\\
\O_{\le n_0} :=\O_{\le n_0} -\O_n,~W:=W-W_n\\
~~~~\textsc{WHILE}~\O_n\ne\emptyset~\textsc{DO}\\
~~~~~~~~~\hbox{Choose any}~o(g_{\ell},g_q)\in\O_{n}\\
~~~~~~~~~\O_n :=\O_n -\{ o(g_{\ell},g_q)\}\\
~~~~~~~~~\OV{o(g_{\ell},g_q)}^{\G}=r\end{array}$\par

$\begin{array}{l} ~~~~~~~~~\textsc{IF}~r\ne 0~\textsc{THEN}\\
~~~~~~~~~~~~~ t':=t'+1,~g_{t'}:=r\\
~~~~~~~~~~~~~\O_{\le n_0} :=\O_{\le n_0}\cup\left\{ o(g_{\ell},g_q) 
~\left |~\begin{array}{l} o(g_{\ell},g_q)\in\left\{\begin{array}{l}
o(g_i,g_{t'}),o(g_{t'},g_i),o(g_{t'},g_{t'}),\\
\hbox{where}~g_i\in\G,~1\le i<t'\end{array}\right\} ,\\
d_{\rm gr}(o(g_{\ell},g_q))\le n_0\end{array}\right.\right\}\\
~~~~~~~~~~~~~\G :=\G\cup\{ g_{t'}\}\\
~~~~~~~~~\textsc{END}\\
~~~~\textsc{END}\\
~~~~\textsc{WHILE}~W_n\ne\emptyset~\textsc{DO}\\
~~~~~~~~~\hbox{Choose any}~f_j\in W_{n}\\
~~~~~~~~~W_n :=W_n -\{ f_j\}\\
~~~~~~~~~\OV{f_j}^{\G}=r\\
~~~~~~~~~\textsc{IF}~r\ne 0~\textsc{THEN}\\
~~~~~~~~~~~~~ t':=t'+1,~g_{t'}:=r\\
~~~~~~~~~~~~~\O_{\le n_0} :=\O_{\le n_0}\cup\left\{ o(g_{\ell},g_q) 
~\left |~\begin{array}{l} o(g_{\ell},g_q)\in\left\{\begin{array}{l}
o(g_i,g_{t'}),o(g_{t'},g_i),o(g_{t'},g_{t'}),\\
\hbox{where}~g_i\in\G,~1\le i<t'\end{array}\right\} ,\\
d_{\rm gr}(o(g_{\ell},g_q))\le n_0\end{array}\right.\right\}\\
~~~~~~~~~~~~~\G :=\G\cup\{ g_{t'}\}\\
~~~~~~~~~\textsc{END}\\
~~~~\textsc{END}\\
\textsc{UNTIL}~\O_{\le n_0}=\emptyset\\
\textsc{END}\end{array}$\par 
\underline{~~~~~~~~~~~~~~~~~~~~~~~~~~~~~~~~~~~~~~~~~~~~~~~~~~~~~~~~~~~~~~~~~~~~~~~~~~~~~~~~~~~~~~~~~~~~~~~~~~~~~~~~~~~~~~~~~~~~ 
~~~~~~} } \vskip 6pt

{\bf Proof} For each fixed $n\le n_0$, by the definition of an
overlap element it is clear that ${\cal O}_n$ is finite. Hence the
algorithm terminates after ${\cal O}_{n_0}$ and $W_{n_0}$ are
exhausted.  Note that  both the WHILE loops append new elements to
$\G$ by taking the nonzero normal remainders on division by $\G$.
With a fixed $n$, by the definition of an overlap element  and the
normality of $g_{t'}$ (mod $\G$), it is straightforward to check
that in both the WHILE loops every nonzero $H\in\{
o(g_i,g_{t'}),o(g_{t'},g_i),o(g_{t'},g_{t'})\}$ has $d_{\rm
gr}(H)>n$. For convenience, let us write $I(n)$ for the ideal
generated by $\G$ which is obtained after $W_n$ is exhausted in the
second WHILE loop. If $n_1$ is the first number after $n$ such that
${\cal O}_{n_1}\ne\emptyset$, and for some $o(g_{\ell},g_q)\in {\cal
O}_{n_1}$, $r =\OV{o(g_{\ell},g_q)}^{\G}\ne 0$ in a certain pass
through the first WHILE loop, then we note that this $r$ is still
contained in $I(n)$. Hence, after ${\cal O}_{n_1}$ is exhausted in
the first WHILE loop, the obtained $\G$ generates $I(n)$ and $\G$ is
an $n_1$-truncated Gr\"obner basis of $I(n)$. Noticing that the
algorithm starts with ${\cal O}=\emptyset$ and $\G =\emptyset$,
inductively it follows from Proposition 2.3 and Corollary 2.4  that
after $W_{n_1}$ is exhausted  in  the second WHILE loop, the
obtained $\G$ is an $n_1$-truncated Gr\"obner basis of $I(n_1)$.
Since $n_0$ is finite and all the generators of $I$ with $d_{\rm
gr}(f_j)\le n_0$ are processed through the second WHILE loop, the
eventually obtained $\G$ is an $n_0$-truncated Gr\"obner basis of
$I$. Finally, the fact that  the degrees of elements in $\G$ are
non-decreasingly ordered follows from the choice of the next $n$ in
the algorithm. \QED \v5

{\bf Remark} Note that in Proposition 2.5 we did not assume that the
subset $F$ is LM-reduced. The reason is that the algorithm starts
with ${\cal O}_{\le n_0}=\emptyset$ and $\G =\emptyset$, while $\G$ 
starts to get its members from the second WHILE loop, and then, the 
new $\G$ obtained after each pass through the WHILE loops is clearly 
LM-reduced. }\v5

Let $I$ be a finitely generated graded ideal of $\KX$. We say that a
homogeneous generating set $F=\{ f_1,\ldots ,f_m\}$ of $I$ is a {\it
minimal homogeneous generating set} if any proper subset of $F$
cannot be a generating set of $I$. We now proceed to show that {\bf
Algorithm 2} presented above can be further modified to compute
minimal homogeneous generating sets for finitely generated graded
ideals of $\KX$. The next proposition and its corollary are
noncommutative analogues of ([KR], Proposition 4.6.1, Corollary
4.6.2).{\parindent=0pt\v5

{\bf 2.6. Proposition} Let $I=\langle F\rangle$ be the graded ideal
of $\KX$ generated by a finite subset of nonzero homogeneous
elements $F=\{ f_1,\ldots ,f_m\}$, where $d_{\rm gr}(f_1)\le d_{\rm
gr}(f_2)\le \cdots\le d_{\rm gr}(f_m)$. Put $I_1=\{ 0\}$,
$I_i=\langle F_i\rangle$, where $F_i=F-\{f_i,\ldots ,f_m\}$,  $2\le
i\le m$. The following statements hold.\par

(i) $F$ is a minimal homogeneous generating set of $I$ if and only
if $f_i\not\in I_i$, $1\le i\le m$.\par

(ii) The set $\OV F=\{ f_k~|~1\le k\le m,~f_k\not\in I_k\}$ is a
minimal homogeneous generating set of $I$.\vskip 6pt

{\bf Proof} (i) If $F$ is a minimal homogeneous generating set of
$I$, then clearly $f_i\not\in I_i$, $1\le i\le m$.}\par

Conversely, suppose $f_i\not\in I_i$, $1\le i\le m$. If $F$ were not
a minimal homogeneous generating set of $I$, then, there is some $i$
such that $I$ is generated by $F'=\{ f_1,\ldots
,f_{i-1},f_{i+1},\ldots ,f_m\}$. Thus, there are $f_j\in F'$ and
nonzero homogeneous elements $h_{jk},h_{j\ell}\in \KX$ such that
$f_i=\sum_{j\ne i}h_{jk}f_jh_{j\ell}$ and $d_{\rm gr}(f_i)=d_{\rm
gr}(h_{jk})+d_{\rm gr}(f_j)+d_{\rm gr}(h_{j\ell})$. Thus $d_{\rm
gr}(f_j)\le d_{\rm gr}(f_i)$ for all $j\ne i$ appeared in the
representation of $f_i$. If $d_{\rm gr}(f_j)<d_{\rm gr}(f_i)$ for
all $j\ne i$, then $f_i\in I_i=\langle F_i\rangle$, which
contradicts the assumption. If $d_{\rm gr}(f_i)=d_{\rm gr}(f_j)$ for
some $j\ne i$, then since $h_{jk}$ and  $h_{j\ell}$ are nonzero
homogeneous elements, we have $h_{jk},h_{j\ell}\in \KX_0-\{
0\}=K^*$. Putting $i'=\max\{i,~j~|~j\ne i,~d_{\rm gr}(f_j)=d_{\rm
gr}(f_i)\}$, we then have $f_{i'}\in I_{i'}=\langle F_{i'}\rangle$,
which again contradicts the assumption. Hence, under the assumption
we conclude that $F$ is a minimal homogeneous generating set of
$I$.\par

(ii) In view of (i), it is sufficient to show that $\OV F$ is a
homogeneous generating set of $I$. Indeed, if $f_i\in F-\OV F$, then
$f_i\in I_i$. By checking $f_{i-1}$ and so on, it follows that
$f_i\in \langle \OV F\rangle$, as desired.\QED{\parindent=0pt\v5

{\bf 2.7. Corollary} Let $F=\{ f_1,\ldots ,f_m\}$ be a minimal
homogeneous generating set of the graded ideal $I$ of $\KX$, where
$d_{\rm gr}(f_1)\le d_{\rm gr}(f_2)\le \cdots\le d_{\rm gr}(f_m)$,
and let $f\in \KX -I$ be a homogeneous element with $ d_{\rm
gr}(f_m)\le d_{\rm gr}(f)$. Then $\widehat{F}= F\cup\{ f\}$ is a
minimal homogeneous generating set of the graded ideal
$\widehat{I}=I+\langle f\rangle$.\par\QED}\v5

Combining the foregoing results,  we  are ready to reach the goal of
this section. {\parindent=0pt\v5

{\bf 2.8. Theorem}  (Compare with ([KR], Theorem 4.6.3)) Let $F =\{ 
f_1,\ldots ,f_m\}$ be a finite set of nonzero homogeneous elements 
of $\KX$ with $d_{\rm gr}(f_1)\le d_{\rm gr}(f_2)\le\cdots \le 
d_{\rm gr}(f_m)=n_0$. Then the following algorithm returns a minimal 
homogeneous generating set $F_{\min}\subseteq F$ for the graded 
ideal $I=\langle F\rangle$; and meanwhile it returns an 
$n_0$-truncated Gr\"obner basis $\G =\{ g_1,\ldots ,g_t\}$ for $I$ 
such that $d_{\rm gr}(g_1)\le d_{\rm gr}(g_2)\le\cdots d_{\rm 
gr}(g_t)$. {\parindent=0pt\vskip 6pt

\underline{\bf Algorithm 3
~~~~~~~~~~~~~~~~~~~~~~~~~~~~~~~~~~~~~~~~~~~~~~~~~~~~~~~~~~~~~~~~~~~~~~~~~~~~~~~~~~~~~~~~~~}\par

$\begin{array}{l} \textsc{INPUT}:~ F= \{ 
f_1,...,f_m\}~\hbox{with}~d_{\rm gr}(f_1)\le d_{\rm 
gr}(f_2)\le\cdots \le 
d_{\rm gr}(f_m)=n_0\\
\textsc{OUTPUT}: ~F_{\min}=\{f_{j_1},\ldots ,f_{j_r}\}\subset F,
~\hbox{a minimal homogeneous generating set of}~I\\
\hskip 2.5truecm\G =\{ g_1,...,g_t\},~\hbox{an}~n_0\hbox{-truncated
Gr\"obner basis
of}~I\\
\textsc{INITIALIZATION}:~ \O_{\le n_0} :=\emptyset ,~W :=F,~\G :=\emptyset ,~t' :=0,~F_{\min}=\emptyset\\
\textsc{LOOP}\\
n:=\min\{d_{\rm gr}(f_i),~d_{\rm gr}(o(g_{\ell},g_q))~|~f_i\in W
,~o(g_{\ell},g_q)\in\O_{\le n_0}\}\\
\O_n:=\{ o(g_{\ell},g_q)\in\O_{\le n_0}~|~d_{\rm gr}(o(g_{\ell},g_q))=n\},~W_n:=\{ f_j\in W~|~d_{\rm gr}(f_j)=n\}\\
\O_{\le n_0} :=\O_{\le n_0} -\O_n,~W:=W-W_n\\
~~~~\textsc{WHILE}~\O_n\ne\emptyset~\textsc{DO}\\
~~~~~~~~~\hbox{Choose any}~o(g_{\ell},g_q)\in\O_{n}\\
~~~~~~~~~\O_n :=\O_n -\{ o(g_{\ell},g_q)\}\\
~~~~~~~~~\OV{o(g_{\ell},g_q)}^{\G}=r\end{array}$\par

$\begin{array}{l}
~~~~~~~~~\textsc{IF}~r\ne 0~\textsc{THEN}\\
~~~~~~~~~~~~~ t':=t'+1,~g_{t'}:=r\\
~~~~~~~~~~~~~\O_{\le n_0} :=\O_{\le n_0}\cup\left\{ o(g_{\ell},g_q) 
~\left |~\begin{array}{l} o(g_{\ell},g_q)\in\left\{\begin{array}{l}
o(g_i,g_{t'}),o(g_{t'},g_i),o(g_{t'},g_{t'}),\\
\hbox{where}~g_i\in\G,~1\le i<t'\end{array}\right\} ,\\
d_{\rm gr}(o(g_{\ell},g_q))\le n_0\end{array}\right.\right\}\\
~~~~~~~~~~~~~\G :=\G\cup\{ g_{t'}\}\\
~~~~~~~~~\textsc{END}\\
\end{array}$\par

$\begin{array}{l}
~~~~\textsc{END}\\
~~~~\textsc{WHILE}~W_n\ne\emptyset~\textsc{DO}\\
~~~~~~~~~\hbox{Choose any}~f_j\in W_{n}\\
~~~~~~~~~W_n :=W_n -\{ f_j\}\\
~~~~~~~~~\OV{f_j}^{\G}=r\\
~~~~~~~~~\textsc{IF}~r\ne 0~\textsc{THEN}\\
~~~~~~~~~~~~~ F_{\min}:=F_{\min}\cup \{f_j\}\\
~~~~~~~~~~~~~ t':=t'+1,~g_{t'}:=r\\
~~~~~~~~~~~~~\O_{\le n_0} :=\O_{\le n_0}\cup\left\{ o(g_{\ell},g_q) 
~\left |~\begin{array}{l} o(g_{\ell},g_q)\in\left\{\begin{array}{l}
o(g_i,g_{t'}),o(g_{t'},g_i),o(g_{t'},g_{t'}),\\
\hbox{where}~g_i\in\G,~1\le i<t'\end{array}\right\} ,\\
d_{\rm gr}(o(g_{\ell},g_q))\le n_0\end{array}\right.\right\}\\
~~~~~~~~~~~~~\G :=\G\cup\{ g_{t'}\}\\
~~~~~~~~~\textsc{END}\\
~~~~\textsc{END}\\
\textsc{UNTIL}~\O_{\le n_0}=\emptyset\\
\textsc{END}\end{array}$\par 
\underline{~~~~~~~~~~~~~~~~~~~~~~~~~~~~~~~~~~~~~~~~~~~~~~~~~~~~~~~~~~~~~~~~~~~~~~~~~~~~~~~~~~~~~~~~~~~~~~~~~~~~~~~~~~~~~~~~~~~~ 
~~~~~~} } \vskip 6pt

{\bf Proof} By  Proposition 2.5 we know that this algorithm
terminates and the eventually obtained $\G$ is an $n_0$-truncated
homogeneous Gr\"obner basis for the ideal $I$, in which the degrees
of elements are ordered non-decreasingly. It remains to prove that
the eventually obtained $F_{\min}$ is a minimal homogeneous
generating set of the ideal $I$.}\par

As in the proof of Proposition 2.5, let us first bear in mind that
for each $n$, in both the WHILE loops every  new appended
$o(g_{\ell},g_{q})$ has $d_{\rm gr}(o(g_{\ell},g_{q}))>n$. Moreover,
for  convenience, let us write $\G (n)$ for the $\G$ obtained after
${\cal  O}_n$ is exhausted in the first WHILE loop, and write
$F_{\min}[n]$, $\G [n]$ respectively  for the $F_{\min}$, $\G$
obtained after $W_n$ is exhausted in the second WHILE loop.  Since
the algorithm starts with ${\cal O}=\emptyset$ and  $\G=\emptyset$,
if, for a fixed $n$, we check carefully how the elements of
$F_{\min}$ are chosen during executing the second WHILE loop, and
how the new elements are appended to $\G$ after each pass through
the first or the second WHILE loop, then it follows from
Proposition 2.3 and Corollary 2.4 that after $W_n$ is exhausted,
the obtained $F_{\min}[n]$ and $\G [n]$ generate the same ideal,
denoted $I(n)$, such that  $\G [n]$ is an $n$-truncated Gr\"obner
basis of $I(n)$.  We now use induction to show that the eventually
obtained  $F_{\min}$ is a minimal homogeneous generating set of the
ideal $I=\langle F\rangle$. If $F_{\min}=\emptyset$, then it is a
minimal generating set of the zero ideal.  To proceed,  we assume
that $F_{\min}[n]$ is a minimal homogeneous generating set for
$I(n)$ after $W_n$ is exhausted in the second WHILE loop. Suppose
that $n_1$ is the first number after $n$ such that ${\cal
O}_{n_1}\ne\emptyset$. We complete the induction proof below by
showing  that $F_{\min}[n_1]$ is a minimal homogeneous generating
set of $I(n_1)$.\par

If in a certain pass through the first WHILE loop,
$r=\OV{o(g_{\ell},g_q)}^{\G}\ne 0$ for some $o(g_{\ell},g_q)\in
{\cal O}_{n_1}$, then we note that $r\in I(n)$. It follows that
after ${\cal O}_{n_1}$ is exhausted in the first WHILE loop, we have
$I(n)=\langle\G (n_1)\rangle$ such that $\G (n_1)$ is an
$n_1$-truncated Gr\"obner basis of $I(n)$. Next, assume that
$W_{n_1}=\{ f_{j_1},\ldots ,f_{j_s}\}\ne\emptyset$ and that the
elements of $W_{n_1}$ are processed in the given order during
executing the second WHILE loop. Since $\G (n_1)$ is an
$n_1$-truncated Gr\"obner basis of $I(n)$, if $f_{j_1}\in W_{n_1}$
is such that $r_1 =\OV {f_{j_1}}^{\G (n_1)}\ne 0$, then $f_{j_1},
r_1\in \KX -I(n)$. By Corollary 2.4, we conclude  that $\G
(n_1)\cup\{r_1\}$ is an $n_1$-truncated Gr\"obner basis for
$I(n)+\langle r_1\rangle$; and by Corollary 2.7, we conclude that
$F_{\min}[n]\cup\{f_{j_1}\}$ is a minimal homogeneous generating set
of $I(n)+\langle r_1\rangle$. Repeating this procedure, if
$f_{j_2}\in W_{n_1}$ is such that $r_2 =\OV{f_{j_2}}^{\G
(n_1)\cup\{r_1\}}\ne 0$, then $f_{j_2}, r_2\in \KX -(I(n)+\langle
r_1\rangle)$.  By Corollary 2.4, we conclude  that $\G
(n_1)\cup\{r_1,r_2\}$ is an $n_1$-truncated Gr\"obner basis for
$I(n)+\langle r_1,r_2\rangle$; and by Corollary 2.7, we conclude
that $F_{\min}[n]\cup\{f_{j_1},f_{j_2}\}$ is a minimal homogeneous
generating set of $I(n)+\langle r_1,r_2\rangle$. Continuing this
procedure until $W_{n_1}$ is exhausted we see that the resulted $\G
[n_1]=\G$ and $F_{\min}[n_1]=F_{\min}$ generate the same module
$I(n_1)$ such that $\G [n_1]$ is an $n_1$-truncated Gr\"obner basis
of $I(n_1)$ and $F_{\min}[n_1]$ is a minimal homogeneous generating
set of $I(n_1)$, as desired. As all elements of $F$ are eventually
processed by the second WHILE loop, we conclude that the finally
obtained $\G$ and $F_{\min}$ have the properties that
$I=\langle\G\rangle$, $\G$ is an $n_0$-truncated Gr\"obner basis of
$I$, and $F_{\min}$ is a minimal homogeneous generating set of
$I$.\QED {\parindent=0pt\v5

{\bf 2.9. Corollary}  Let $F =\{ f_1,\ldots ,f_m\}$ be a finite set
of nonzero homogeneous elements of $\KX$ with $d_{\rm gr}(f_1)=
d_{\rm gr}(f_2)=\cdots = d_{\rm gr}(f_m)=n_0$. \par

(i) If  $\LM (f_i)\ne\LM (f_j)$ for all $i\ne j$, then $F$ is a 
minimal homogeneous generating set of the ideal $I=\langle 
F\rangle$, and meanwhile $F$ is an $n_0$-truncated Gr\"obner basis 
for $I$.\par

(ii) If $F$ is a minimal Gr\"obner basis of the ideal $I=\langle
F\rangle$, then  $F$ is a minimal homogeneous generating set of
$I$.\vskip 6pt

{\bf Proof} By the assumption, it follows from the second WHILE loop
of {\bf Algorithm 3} that $F_{\min}=F$.} \v5

\section*{3. Computation of Minimal Standard Bases}
Let $\KX =K\langle X_1,\ldots ,X_n\rangle$ be the free $K$-algebra
generated by $X=\{ X_1,\ldots ,X_n\}$ and $\B$ the standard
$K$-basis of $\KX$. Fix a weighted $\NZ$-gradation $\KX
=\oplus_{q\in\NZ}\KX_q$ for $\KX$ by assigning each $X_i$ a positive
degree $d_{\rm gr}(X_i)=m_i$, $1\le i\le n$. Recall that a {\it
graded monomial ordering} on $\B$ is a monomial ordering $\prec$ on
$\B$ satisfying
$$u,v\in\B~\hbox{and}~u\prec v~\hbox{implies}~d_{\rm gr}(u)\le d_{\rm gr}(v).$$
A graded monomial ordering is usually denoted by $\prec_{gr}$. The
most well-known graded monomial ordering on $\B$ is the graded
lexicographic ordering $\prec_{grlex}$.\par

In this section,  we show that if an {\it ungraded ideal} $I$ of
$\KX$ has a finite Gr\"obner basis $\G$ with respect to a given
graded monomial ordering $\prec_{gr}$, then a minimal standard basis
for $I$  can be computed via computing a minimal homogeneous
generating set of the associated graded ideal $\langle\LH
(I)\rangle$ of $I$ (see the definitions below). Concerning the
notion of a standard basis for the ideal $I$, we have a remark given
after Proposition 3.2 below.   All notions, notations, and
conventions used before are maintained.\v5

Let $f=f_0+f_1+\cdots +f_q\in \KX$ with $f_i\in \KX_i$ and $f_q\ne
0$, and let $\LH (f)$ denote the {\it leading homogeneous element}
of $f$, i.e., $\LH (f)=f_q$. Then every ideal $I$ of $\KX$ has the
associated graded ideal $\langle\LH (I )\rangle$ generated by the
set of leading homogeneous elements $\LH (I)=\{ \LH (f)~|~f\in
I\}$.{\parindent=0pt\v5

{\bf 3.1. Definition} Let $I$ be an arbitrary ideal of $\KX$. A
subset $G$ of $I$ is said to be a {\it standard basis} for $I$, if
$\langle \LH (I)\rangle =\langle \LH (G)\rangle$.\v5

{\bf 3.2. Proposition} With respect to the fixed weighted
$\NZ$-graded $K$-algebra structure $\KX =\oplus_{q\in\NZ}\KX_q$, let
$\KX$ be equipped with the $\NZ$-grading filtration $F\KX =\{
F_q\KX\}_{q\in\NZ}$, where for each $q\in\NZ$, $F_q\KX =\oplus_{k\le
q}\KX_k$, and let $I$ be an arbitrary ideal of $\KX$. For a subset
$G$ of $I$, the following  statements are equivalent.\par

(i) $G$ is a standard basis of $I$;\par

(ii) Every nonzero element $f\in I$ has a representation
$$\begin{array}{rcl} f&=&\sum_{i,j} \lambda_{ij}u_{ij}g_jv_{ij},~\lambda_{ij}\in K,~u_{ij},v_{ij}\in\B,\\
&{~}&\hbox{satisfying}~d_{\rm gr}(\LH (u_{ij}g_jv_{ij}))\le d_{\rm
gr}(\LH (f))~\hbox{for all}~\lambda_{ij}\ne 0;\end{array}$$\par

(iii) Let $d_{\rm gr}(g_j)=q_j$, $g_j\in G$. Considering the induced
filtration $FI=\{F_qI\}_{q\in\NZ}$ of $I$ with $F_qI=I\cap F_q\KX$,
we have
$$ F_qI=\sum_{g_j\in G}\left (\sum_{k_i+q_j+k_j\le q}F_{k_i}\KX g_jF_{k_j}\KX\right ),\quad q\in\NZ.$$\par

{\bf Proof} This is referred to the proof of ([LWZ], Lemma
2.2.3).\QED}\v5

By Proposition 3.2 it is clear that every standard basis $G$ of $I$
is certainly a generating set of $I$. By Definition 3.1 it is also
clear that if $I$ is a graded ideal of $\KX$, then any homogeneous
generating set $G$ of $I$ is trivially a standard basis of $I$.
Nevertheless, we shall continue our discussion below for arbitrary
ideals. Moreover, we specify the following{\parindent=0pt\v5

{\bf Remark} As one may see from the literature on computational
commutative algebra (e.g. see [KR]), if $A=K[x_1,\ldots ,x_n]$ is
the commutative polynomial $K$-algebra in $n$ variables, then a
standard basis for an ideal $I$ of $A$ is nothing but the well-known
Macaulay basis. While in the noncommutative case,  for two-sided
ideals of a $\Gamma$-filtered algebra $A$, where $\Gamma$ is an
ordered semigroup with respect to a well-ordering, standard bases
were introduced in [Gol] by using the induced filtration and the
associated graded ideals. When a weighted $\NZ$-gradation $\KX
=\oplus_{q\in\NZ}\KX_q$ is fixed for the free algebra $\KX =K\langle
X_1,\ldots ,X_n\rangle$, and furthermore $\KX$ is equipped with the
$\NZ$-grading filtration $F\KX =\{ F_q\KX\}_{q\in\NZ}$, where
$F_q\KX =\oplus_{k\le q}\KX_k$, the definition of a standard basis
in the sense of [Gol] is then turned out to be Definition 3.1 above
by Proposition 3.2. In this case, if $G$ is a standard basis of an
ideal $I$ in $\KX$ and if the quotient algebra $A=\KX /I$ is
equipped with the  filtration $FA$ induced by $F\KX$, then the
$\NZ$-filtered algebra $A$ has the associated graded algebra
$G(A)\cong \KX /\langle\LH (G)\rangle$. So, among other
applications, the structure of standard bases for ideals of $\KX$
plays an important role in the study of general PBW theory and the
study of homogeneous and inhomogeneous Koszul algebras. On this
aspect one may refer to ([Li2], Chapter 4) for more details. }\v5

Actually as in the commutative case with a Macaulay basis, we have
the following{\parindent=0pt\v5

{\bf 3.3. Proposition} Let $\prec_{gr}$ be a graded monomial
ordering on $\B$ as defined in the beginning of this section, and
let $I$ be an ideal of $\KX$. If $\G$ is a Gr\"obner basis for $I$
with respect to $\prec_{gr}$, then $\G$ is a standard basis for $I$
in the sense of Definition 3.1, i.e., $\langle\LH (I)\rangle
=\langle\LH (\G )\rangle$.
\par\QED}\v5

Let $I$ be an ideal of $\KX$. If any proper subset of a standard
basis $G$ of $I$ cannot be a standard basis for $I$, then $G$ is
called a {\it minimal standard basis}. By Definition 3.1 it is clear
that a subset $G$ of $I$ is a minimal standard basis for $I$ if and
only if $\LH (G)$ is a minimal homogeneous generating set of the
graded ideal $\langle\LH (I)\rangle$. Thus, as with minimal
homogeneous generating sets for graded ideals, minimal standard
bases have the following properties: {\parindent=1.3truecm\par

\item{(1)} any two minimal standard bases of $I$ have the same number of
generators; and

\item{(2)} any two minimal standard bases of $I$ contain the same
number of leading homogeneous elements of degree $n$ for all
$n\in\NZ$.\v5}

Now, it follows from Proposition 3.3 and Theorem 2.8 that we are
able to give the main result of this section.{\parindent=0pt\v5

{\bf 3.4. Theorem} Let $\prec_{gr}$ be a graded monomial ordering on
$\B$ as defined in the beginning of this section, and let $I$ be an
ideal of $\KX$. If $\G=\{ g_1,\ldots ,g_m\}$ is a finite Gr\"obner
basis for $I$ with respect to $\prec_{gr}$, then a minimal standard
basis of $I$ can be computed by following the steps below:} \par

{\bf Step 1.} With the initial input data $F=\{\LH (g_1),\ldots ,\LH
(g_m)\}$, run {\bf Algorithm 3} to compute a minimal homogeneous
generating set $F_{\min}$ for the graded ideal $\langle\LH
(I)\rangle$, say $F_{\min}=\{ \LH (g_{j_1}),\ldots ,\LH
(g_{j_s}\}$.\par

{\bf Step 2.} Write down  $G=\{ g_{j_1},\ldots ,g_{j_s}\}$, that is
a minimal standard basis of $I$.\par\QED\v5

It follows from Corollary 2.9 and Theorem 3.4 that we have also the
following{\parindent=0pt\v5

{\bf 3.5. Corollary}  Let $I$ be an ideal of $\KX$ and let $\G =\{
g_1,\ldots ,g_m\}$ be a finite  Gr\"obner basis of $I$ with respect
to a graded monomial ordering $\prec_{gr}$ on $\B$.  If $\G$ is a
minimal Gr\"obner basis and $d_{\rm gr}(\LH (g_1))= d_{\rm gr}(\LH
(g_2))=\cdots = d_{\rm gr}(\LH (g_m))=n_0$, then $\G$ is a minimal
standard basis for $I$.\par\QED}\v5

Finally, in the light of Gr\"obner basis theory for path algebras
(i.e. quiver algebras) [Gr], we remark that the results obtained in
this paper hold true for path algebras  defined by finite directed
graphs.

\v5 \centerline{References}\par \parindent=1truecm\par

\item{[Ber]} G. Bergman, The diamond lemma for ring theory, {\it
Adv. Math}., 29(1978), 178--218.

\item{[CDNR]} A. Capani, G. De Dominicis, G. Niesi, and L. Robbiano,
Computing minimal finite free resolutions. {\it Journal of Pure and
Applied Algebra}, (117\& 118)(1997), 105 -- 117.

\item{[Coc]} CoCoATeam,
{\it CoCoA: a system for doing Computations in Commutative Algebra}.
Available at http://cocoa.dima.unige.it

\item{[Gol]} E. S. Golod, Standard bases and homology, in: {\it Some
Current Trends in Algebra}, (Varna, 1986). Lecture Notes in
Mathematics, Vol. 1352, Springer-Verlag, 1988, 88-95.

\item{[Gr]} E. L. Green, Noncommutative Gr¡§obner bases and projective
resolutions, in: {\it Proceedings of the Euroconference
Computational Methods for Representations of Groups and Algebras},
Essen, 1997, (Michler, Schneider, eds). Progress in Mathematics,
Vol. 173, Basel, Birkha¡§user Verlag, 1999, 29--60.

\item{[KR]} M. Kreuzer, L. Robbiano, {\it Computational Commutative Algebra 2}. Springer, 2005.

\item{[Li1]} H. Li, On monoid graded local rings. {\it Journal of Pure and Applied Algebra}, 216(2012), 2697 -- 2708.

\item{[Li2]} H. Li, {\it Gr\"obner Bases in Ring Theory}. World Scientific Publishing Co., 2011.

\item{[LWZ]} H. Li, Y. Wu and J. Zhang, Two applications of noncommutative Gr\"obner bases. {\it
 Annali dell'Universit\'a di Ferrara. Sezione 7: Scienze matematiche}, 45(1)(1999), 1-24.

\item{[Mor]} T. Mora, An introduction to commutative and
noncommutative Gr\"obner bases. {\it Theoretic Computer Science},
134(1994), 131--173.

\end{document}